\documentclass[pdflatex,sn-mathphys-num]{sn-jnl}

\usepackage{graphicx}%
\usepackage{multirow}%
\usepackage{amsmath,amssymb,amsfonts}%
\usepackage{amsthm}%
\usepackage{mathrsfs}%
\usepackage[title]{appendix}%
\usepackage{xcolor}%
\usepackage{textcomp}%
\usepackage{manyfoot}%
\usepackage{booktabs}%
\usepackage{algorithm}%
\usepackage{algorithmicx}%
\usepackage{algpseudocode}%
\usepackage{listings}%

\theoremstyle{thmstyleone}%
\newtheorem{theorem}{Theorem}
\newtheorem{proposition}[theorem]{Proposition}%

\theoremstyle{thmstyletwo}%
\newtheorem{remark}{Remark}%

\theoremstyle{thmstylethree}%

\newcommand{\R}{\mathbb{R}}
\DeclareMathOperator*{\argmin}{arg\,min}
\DeclareMathOperator*{\argmax}{arg\,max}
\newcommand{\norm}[1]{\left\|#1\right\|}
\hypersetup{hypertexnames=false}

\begin{document}

\title[A Non-ergodic Convergence Rate of ADMM for Cosparse Signal Reconstruction]{A Non-ergodic Convergence Rate of ADMM for Cosparse Signal Reconstruction}


\author[1]{\fnm{Ting} \sur{Zhang}}\email{zhangting314314@126.com}

\author*[2]{\fnm{Zisheng} \sur{Liu}}\email{liuzisheng0710@163.com}

\author[3]{\fnm{Yuan} \sur{Liu}}\email{584563738@qq.com}

\author[3]{\fnm{Yonglong} \sur{Ma}}\email{1295906656@qq.com}

\affil[1]{\orgdiv{School of Mathematics and Information Science}, \orgname{Henan University of Economics and Law}, \orgaddress{ \city{Zhengzhou}, \postcode{450046}, \country{China}}}

\affil*[2]{\orgdiv{School of Computational Economics}, \orgname{Henan University of Economics and Law}, \orgaddress{\city{Zhengzhou}, \postcode{450046}, \country{China}}}

\affil[3]{\orgdiv{School of Statistics}, \orgname{Henan University of Economics and Law}, \orgaddress{\city{Zhengzhou}, \postcode{450046}, \country{China}}}


\abstract{This paper investigates the non-ergodic convergence rate of the alternating direction
method of multipliers (ADMM) for cosparse signal reconstruction in a uniformly tight frame.
While the ergodic $O(1/t)$ convergence rate of ADMM for the cosparse optimization problem
has been established in previous work, the non-ergodic convergence behavior remains
unexplored. In this paper, we first introduce auxiliary variables and formulate the cosparse
optimization problem as a variational inequality. Then, based on the strict contraction
property of the ADMM iterative sequence, we prove that the sequence is monotonically
non-increasing with respect to a weighted norm. Using this monotonicity property together
with the strict contraction theorem, we establish a worst-case $O(1/t)$ convergence rate
in the non-ergodic sense, which guarantees that the current iterate itself converges to the
optimal solution set at this rate. This result complements the existing ergodic convergence
analysis and provides a stronger theoretical guarantee for practical ADMM implementations.
Finally, numerical experiments verify the effectiveness of the ADMM method and the
correctness of the non-ergodic convergence theory.}

\keywords{cosparse signal reconstruction, ADMM, uniformly tight frame,
non-ergodic convergence rate, variational inequality}


\pacs[MSC Classification]{90C25, 90C33}

\maketitle

\section{Introduction}

As an alternative to mature sparse representation models, cosparse analysis
models~\cite{nam2011,nam2013,giryes2014} have attracted widespread attention in signal
processing and compressed sensing. Consider the linear system
\begin{equation}\label{eq:1}
    y = Mx + v,
\end{equation}
where $y \in \R^m$ is the observation vector, $M \in \R^{m\times d}$ is the measurement
matrix, and $v$ is the noise vector. Let the redundant analysis operator
$D \in \R^{n\times d}$ ($n \ge d$) be a uniformly tight frame. A signal $x$ is
$\ell$-cosparse under $D$ if $\ell = n - \norm{Dx}_0$ is sufficiently large.
The cosparse signal reconstruction problem can be formulated as the following
optimization problem:
\begin{equation}\label{eq:2}
    \min_x \norm{Dx}_0 \quad \text{s.t.} \quad y = Mx + v.
\end{equation}
Problem~\eqref{eq:2} is NP-hard due to the $\ell_0$ norm.
Typical methods include variants of the subgradient algorithm (CIF \cite{wang2016}), greedy pursuit (such as GAP~\cite{nam2013} and ReGAP~\cite{liu2017,li2015}), iterative methods inspired by CoSaMP and SP, and iterative hard threshold variants \cite{giryes2014,tirer2020,kotzagiannidis2022}. These algorithms belong to the greedy class and aim to find the cosparse representation of the signal while dealing with noise and uncertainty in the measurement process.

Another type of method is the convex relaxation method. Based on the fact that a
non-convex $\ell_0$-norm can be approximated by a convex $\ell_1$-norm, the following
convex relaxation problem is obtained:
\begin{equation}\label{eq:3}
    \min_x \norm{Dx}_1 \quad \text{s.t.} \quad \norm{y - Mx}_2 \le \epsilon,
\end{equation}
where $\epsilon$ is the upper bound of the noise level $\norm{v}_2$. The method for
solving \eqref{eq:3} is called Analysis Basis Pursuit (ABP~\cite{zhao2014}).

Inspired by \cite{he2012,he2018,he2022b}, this paper investigates the non-ergodic convergence of ADMM for solving the cosparse
signal reconstruction problem within the variational inequality (VI) framework.
The main contributions include:
(i) establishing the monotonic non-increasing property of the ADMM sequence for the
cosparse problem under a weighted norm;
(ii) deriving a worst-case $O(1/t)$ convergence rate in the non-ergodic sense;
(iii) the results do not depend on the
RIP~\cite{candes2004,candes2006,baraniuk2008,foucart2013,qaisar2013} or
$D$-RIP~\cite{zhao2014,giryes2014,candes2011,liu2017,liu2026} conditions;
(iv) since $D$ is a uniformly tight frame ($D^T D = I$), the conclusions are not
constrained by symmetric positive semidefiniteness assumptions.
To apply ADMM, we reformulate \eqref{eq:3} as the following equivalent constrained
optimization problem:
\begin{equation}\label{eq:4}
    \min_{x,z} \norm{z}_1 + \frac{\alpha}{2}\norm{y - Mx}_2^2
    \quad \text{s.t.} \quad Dx - z = 0,
\end{equation}
where an auxiliary variable $z \in \R^n$ is introduced, and $\alpha > 0$ is the
penalty parameter.

The paper is organized as follows:
Section~\ref{sec:prelim} presents the VI characterization;
Section~\ref{sec:ergodic} briefly reviews the ergodic convergence results from~\cite{liu2024};
Section~\ref{sec:nonergodic} establishes the main results: the monotonic non-increasing
property and the worst-case $O(1/t)$ non-ergodic convergence rate;
Section~\ref{sec:experiments} presents numerical experiments;
Section~\ref{sec:conclusions} concludes.

\section{Preliminaries: Variational Inequality Characterization}
\label{sec:prelim}

The Lagrangian function of problem~\eqref{eq:4} is:
\begin{equation}\label{eq:5}
    L(z,x,\lambda) = \norm{z}_1 + \frac{\alpha}{2}\norm{y - Mx}_2^2
    - \lambda^T(Dx - z),
\end{equation}
where $\lambda \in \R^n$ is the Lagrange multiplier. We assume $x \in X$, $z \in Z$,
and $\lambda \in \R^n$, where $X \subset \R^d$ and $Z \subset \R^n$ are closed convex sets.
A saddle point $(z^*, x^*, \lambda^*) \in \Omega := Z \times X \times \R^n$ of
$L(z,x,\lambda)$ satisfies:
\begin{equation}\label{eq:6}
    L(z^*, x^*, \lambda) \le L(z^*, x^*, \lambda^*) \le L(z, x, \lambda^*).
\end{equation}
Equivalently, the saddle point can be characterized by the system:
\begin{equation}\label{eq:7}
\begin{cases}
    z^* = \argmin_{z\in Z} L(z, x^*, \lambda^*),\\
    x^* = \argmin_{x\in X} L(z^*, x, \lambda^*),\\
    \lambda^* = \argmax_{\lambda\in\R^n} L(z^*, x^*, \lambda).
\end{cases}
\end{equation}
The following proposition provides the optimality condition via variational inequality.

\begin{proposition}[\cite{he2016}]\label{prop:vi}
Let $X \subset \R^d$ be a closed convex set, and $\theta(x): \R^d \to \R$ be a convex
function. Furthermore, let $f(x)$ be differentiable on $X$. Assume that the solution
set of $\min\{\theta(x) + f(x) \mid x \in X\}$ is nonempty. Then,
\[
x^* = \argmin_{x\in X}\{\theta(x) + f(x)\}
\]
if and only if
\begin{equation}\label{eq:8}
    x, x^* \in X,\quad \theta(x) - \theta(x^*) + (x - x^*)^T \nabla f(x^*) \ge 0.
\end{equation}
\end{proposition}

Let $\theta_1(z) = \norm{z}_1$ and $\theta_2(x) = \frac{\alpha}{2}\norm{y - Mx}_2^2$.
According to Proposition~\ref{prop:vi}, a saddle point $(z^*, x^*, \lambda^*)$ of the
Lagrangian~\eqref{eq:5} can be characterized by a solution of the following variational
inequality: find $\omega^* \in \Omega$ such that
\begin{equation}\label{eq:9}
    \theta(u) - \theta(u^*) + (\omega - \omega^*)^T F(\omega^*) \ge 0,
    \quad \forall \omega \in \Omega,
\end{equation}
where
\begin{equation}\label{eq:10}
    \theta(u) = \theta_1(z) + \theta_2(x),\quad \Omega = Z \times X \times \R^n,
\end{equation}
and
\begin{equation}\label{eq:11}
    \omega = \begin{pmatrix} z \\ x \\ \lambda \end{pmatrix},\quad
    u = \begin{pmatrix} z \\ x \end{pmatrix},\quad
    F(\omega) = \begin{pmatrix} \lambda \\ -D^T\lambda \\ Dx - z \end{pmatrix}.
\end{equation}
Since $F$ is an affine operator, we have
\begin{equation}\label{eq:12}
    F(\omega) = \begin{pmatrix}
        0 & 0 & I \\
        0 & 0 & -D^T \\
        -I & D & 0
    \end{pmatrix}
    \begin{pmatrix} z \\ x \\ \lambda \end{pmatrix}.
\end{equation}
According to the antisymmetry of the affine matrix,
\begin{equation}\label{eq:13}
    (\omega - \bar{\omega})^T [F(\omega) - F(\bar{\omega})] \equiv 0,
    \quad \forall \omega, \bar{\omega} \in \Omega.
\end{equation}

\section{Ergodic Convergence Results (Previous Work)}
\label{sec:ergodic}

This section briefly reviews the ergodic convergence results established in~\cite{liu2024}.
These results serve as the foundation for the non-ergodic analysis in
Section~\ref{sec:nonergodic}.

\subsection{Review of Known Results}

The augmented Lagrangian function of problem~\eqref{eq:4} is:
\begin{equation}\label{eq:14}
    L_{\beta}(z,x,\lambda) = \norm{z}_1 + \frac{\alpha}{2}\norm{y - Mx}_2^2
    - \lambda^T(Dx - z) + \frac{\beta}{2}\norm{Dx - z}_2^2,
\end{equation}
where $\beta > 0$ is the penalty parameter. The ADMM iterative scheme is:
\begin{equation}\label{eq:15}
\begin{cases}
    z^{k+1} = \argmin_{z\in Z} L_{\beta}(z, x^k, \lambda^k),\\[4pt]
    x^{k+1} = \argmin_{x\in X} L_{\beta}(z^{k+1}, x, \lambda^k),\\[4pt]
    \lambda^{k+1} = \lambda^k - \beta(Dx^{k+1} - z^{k+1}).
\end{cases}
\end{equation}
The corresponding variational inequality of~\eqref{eq:15} is:
\begin{equation}\label{eq:16}
\begin{cases}
    \theta_1(z) - \theta_1(z^{k+1}) + (z - z^{k+1})^T
        [\lambda^k - \beta(Dx^k - z^{k+1})] \ge 0,\\[4pt]
    \theta_2(x) - \theta_2(x^{k+1}) + (x - x^{k+1})^T
        [-D^T\lambda^k + \beta D^T(Dx^{k+1} - z^{k+1})] \ge 0,\\[4pt]
    (\lambda - \lambda^{k+1})^T
        [(Dx^{k+1} - z^{k+1}) + \tfrac{1}{\beta}(\lambda^{k+1} - \lambda^k)] \ge 0.
\end{cases}
\end{equation}
In~\cite{liu2024}, three technical lemmas (Lemma~3.1--3.3) were established to demonstrate
the strict contraction property of the ADMM sequence. Based on these lemmas, the following
strict contraction theorem was proved.

\begin{theorem}[Strict Contraction, \cite{liu2024}]\label{thm:strict}
Let the sequence $\{\omega^k\}$ be generated by~\eqref{eq:15}. Then, for any saddle
point $(z^*, x^*, \lambda^*)$, we have
\begin{equation}\label{eq:20}
    \norm{v^{k+1} - v^*}_H^2 \le
    \norm{v^k - v^*}_H^2 - \norm{v^k - v^{k+1}}_H^2,
\end{equation}
where
\begin{equation}\label{eq:21}
    v = \begin{pmatrix} \lambda \\ x \end{pmatrix},\quad
    H = \begin{pmatrix} \frac{1}{\beta} I_m & 0 \\ 0 & \beta I_d \end{pmatrix},
\end{equation}
and $V^* = \{(\lambda^*, x^*) \mid (z^*, x^*, \lambda^*) \in \Omega\}$.
\end{theorem}

Since $D$ is a uniformly tight frame ($D^T D = I$), the matrix $H$ is symmetric positive
definite. From Theorem~\ref{thm:strict}, the sequence $\{v^k\}$ is strictly contractive
with respect to the solution set $V^*$, and thus converges to a point in $V^*$.
Summing both sides of~\eqref{eq:20} yields:
\begin{equation}\label{eq:22}
    \sum_{k=0}^{\infty} \norm{v^k - v^{k+1}}_H^2 \le \norm{v^0 - v^*}_H^2,
\end{equation}
which implies $\lim_{k\to\infty} \norm{v^k - v^{k+1}}_H^2 = 0$.

\subsection{Ergodic Convergence Rate}

In~\cite{liu2024}, the worst-case $O(1/t)$ convergence rate in the ergodic sense was
also established as follows.

\begin{theorem}[\cite{liu2024}]\label{thm:ergodic}
Let $\{\omega^k\}$ be generated by~\eqref{eq:15}. Then for any integer $t \ge 0$,
\begin{equation}\label{eq:23}
    \theta(\bar{u}^t) - \theta(u) + (\bar{\omega}^t - \omega)^T F(\omega)
    \le \frac{1}{2(t+1)}
    \left[\frac{1}{\beta}\norm{\lambda^0 - \lambda}_2^2
    + \beta\norm{Dx^0 - z}_2^2\right],
\end{equation}
where
\begin{equation}\label{eq:24}
    \bar{\omega}^t = \frac{1}{t+1}\sum_{k=0}^t \omega^{k+1},
\end{equation}
and $\omega \in \Omega$. This establishes the worst-case $O(1/t)$ convergence rate of
ADMM for the cosparse signal reconstruction problem in the ergodic sense.
\end{theorem}

\section{Convergence Rate in the Non-ergodic Sense}
\label{sec:nonergodic}

In this section, we establish the main result: the worst-case $O(1/t)$ convergence rate
of ADMM for cosparse signal reconstruction in the non-ergodic sense. Unlike the ergodic
rate (which concerns the averaged sequence $\bar{\omega}^t$), the non-ergodic rate
guarantees that the current iterate $\omega^{k+1}$ converges to the optimal solution set
at rate $O(1/t)$. This is practically more relevant because algorithms output the current
iterate, not an average. To derive this rate, we first prove that
$\norm{v^k - v^*}_H^2$ is monotonically non-increasing.

\begin{theorem}\label{thm:monotone}
Assume that the sequence $\{v^k\}$ is generated by~\eqref{eq:15}. Then,
\begin{equation}\label{eq:25}
    \norm{v^{k+1} - v^*}_H^2 \le \norm{v^k - v^*}_H^2,
\end{equation}
where $H$ is as defined in~\eqref{eq:21}.
\end{theorem}

\begin{proof}
According to Lemma~3.1 of~\cite{liu2024}, we have
\begin{equation}\label{eq:26}
\begin{aligned}
&\theta(u) - \theta(u^{k+1}) + (\omega - \omega^{k+1})^T F(\omega)\\
\geq& (z - z^{k+1})^T \beta (Dx^k - Dx^{k+1}) + \frac{1}{\beta}(\lambda - \lambda^{k+1})^T(\lambda^k - \lambda^{k+1}), \\
=& \beta[(z - z^{k+1}) + (Dx^{k+1} - Dx)]^T (Dx^k - Dx^{k+1}) \\
+&  \beta(Dx - Dx^{k+1})^T (Dx^k - Dx^{k+1}) + \frac{1}{\beta}(\lambda - \lambda^{k+1})^T(\lambda^k - \lambda^{k+1}),~\forall\omega\in\Omega.
\end{aligned}
\end{equation}
According to equation~\eqref{eq:21} , equation~\eqref{eq:26} can be rewritten as follows:
\begin{equation}\label{e63}
\begin{aligned}
&\theta(u) - \theta(u^{k+1}) + (\omega - \omega^{k+1})^T F(\omega)\\
\geq&  \beta[(z - z^{k+1}) + (Dx^{k+1} - Dx)]^T (Dx^k - Dx^{k+1})  \\
+&  (v - v^{k+1})^T H (v^k - v^{k+1}),~\forall\omega\in\Omega.
\end{aligned}
\end{equation}
Let  $\omega = \omega^k$, according to the following formula
\begin{eqnarray}\label{e64}
\beta(Dx^{k+1} - z^{k+1}) = (\lambda^k - \lambda^{k+1}),~\textrm{and}~ \beta(Dx^{k} - z^{k}) = (\lambda^{k-1} - \lambda^{k}),
\end{eqnarray}
Then equation (\ref{e63})  can be rewritten as
\begin{equation}\label{e65}
\begin{aligned}
&\theta(u^k) - \theta(u^{k+1}) + (\omega^k - \omega^{k+1})^T F(\omega^k)\\
\geq&  [(\lambda^k - \lambda^{k+1}) - (\lambda^{k-1} - \lambda^k)]^T (Dx^k - Dx^{k+1})  \\
+ & (v^k - v^{k+1})^T H (v^k - v^{k+1}),~\forall\omega\in\Omega.
\end{aligned}
\end{equation}
Similarly, if we let  $k = k-1$ in equation (\ref{e63}), then equation (\ref{e63}) can be rewritten as
\begin{equation}\label{e66}
\begin{aligned}
&\theta(u) - \theta(u^{k}) + (\omega - \omega^{k})^T F(\omega)\\
\geq&  \beta[(z - z^{k}) + (Dx^{k} - Dx)]^T (Dx^{k-1} - Dx^{k}) \\
+& (v - v^{k})^T H (v^{k-1} - v^{k}),~\forall\omega\in\Omega.
\end{aligned}
\end{equation}
Let  $\omega = \omega^{k+1}$, and using equation  (\ref{e64}), then equation  (\ref{e66}) can be rewritten as
\begin{equation}\label{e67}
\begin{aligned}
&\theta(u^{k+1}) - \theta(u^{k}) + (\omega^{k+1} - \omega^{k})^T F(\omega^{k+1})\\
\geq& [(\lambda^{k-1} - \lambda^{k}) - (\lambda^{k} - \lambda^{k+1})]^T (Dx^{k-1} - Dx^{k})  \\
+& (v^{k+1} - v^{k})^T H (v^{k-1} - v^{k}),~\forall\omega\in\Omega.
\end{aligned}
\end{equation}
Combine inequalities \eqref{e65} and \eqref{e67}. By utilizing the identity $(\omega^k - \omega^{k+1})^T F(\omega^k - \omega^{k+1}) = 0$, we have
\begin{equation}\label{e68}
\begin{aligned}
0\geq& [(\lambda^k - \lambda^{k+1}) - (\lambda^{k-1} - \lambda^k)]^T (Dx^k - Dx^{k+1})  + (v^k - v^{k+1})^T H (v^k - v^{k+1})\\
 +& [(\lambda^{k-1} - \lambda^{k}) - (\lambda^{k} - \lambda^{k+1})]^T (Dx^{k-1} - Dx^{k})  + (v^{k+1} - v^{k})^T H (v^{k-1} - v^{k}).
\end{aligned}
\end{equation}
This means
\begin{equation}\label{e69}
\begin{aligned}
&(v^k - v^{k+1})^T H [(v^{k-1} - v^{k}) - (v^k - v^{k+1})]\\
\geq& [(\lambda^{k-1} - \lambda^{k}) - (\lambda^{k} - \lambda^{k+1})]^T D [(x^{k-1} - x^{k}) - (x^{k} - x^{k + 1})].
\end{aligned}
\end{equation}
Let
\begin{eqnarray}\label{e70}
a = (v^{k-1} - v^{k}) ~\textrm{and} ~ b = (v^k - v^{k+1}),
\end{eqnarray}
based on the inequality
\begin{eqnarray}\label{e71}
\|a\|_H^2 - \|b\|_H^2 = \|a - b\|_H^2 + 2b^T H(a - b),
\end{eqnarray}
we have
\begin{equation}\label{e72}
\begin{aligned}
&\|v^{k-1} - v^{k}\|_H^2 - \|v^k - v^{k+1}\|_H^2 \\
=& \|(v^{k-1} - v^{k}) - (v^k - v^{k+1})\|_H^2 + 2(v^k - v^{k+1})^T H[(v^{k-1} - v^{k}) - (v^k - v^{k+1})].
\end{aligned}
\end{equation}
Using equation \eqref{e69}, equation \eqref{e72} can be rewritten as follows:
\begin{equation}\label{e73}
\begin{aligned}
&\|v^{k-1} - v^{k}\|_H^2 - \|v^k - v^{k+1}\|_H^2 \\
\geq& \|(v^{k-1} - v^{k}) - (v^k - v^{k+1})\|_H^2 \\
+& 2[(\lambda^{k-1} - \lambda^{k}) - (\lambda^{k} - \lambda^{k+1})]^T D [(x^{k-1} - x^{k}) - (x^{k} - x^{k + 1})].
\end{aligned}
\end{equation}
Based on the symmetric positive semi-definite structure of matrix $H$   and the sum of squares formula, the right side of equation \eqref{e73} equals
\begin{eqnarray}\label{e74}
\|\sqrt{\beta}D[(x^{k-1} - x^{k}) - (x^{k} - x^{k + 1})] + \frac{1}{\sqrt{\beta}}[(\lambda^{k-1} - \lambda^{k}) - (\lambda^{k} - \lambda^{k+1})]\|^2.
\end{eqnarray}
Since equation \eqref{e74} is nonnegative, the theorem is proved.
\end{proof}

\begin{remark}
Theorem~\ref{thm:monotone} establishes a stronger property than the strict contraction
result of Theorem~\ref{thm:strict}. While Theorem~\ref{thm:strict} shows that the
distance to the solution set decreases at each iteration, Theorem~\ref{thm:monotone}
further shows that the sequence $\{\norm{v^k - v^*}_H^2\}$ is monotonic, which is the
key enabling factor for the non-ergodic convergence rate.
\end{remark}

Combining Theorems~\ref{thm:strict} and~\ref{thm:monotone}, we now prove the worst-case
$O(1/t)$ convergence rate in the non-ergodic sense.

\begin{theorem}\label{thm:nonergodic}
Assume that the sequence $\{v^k\}$ is generated by~\eqref{eq:15}. Then,
\begin{equation}\label{eq:39}
    \norm{v^t - v^{t+1}}_H^2
    \le \frac{1}{t+1}\norm{v^0 - v^*}_H^2,
\end{equation}
where $H$ and $v^*$ are given by Theorem~\ref{thm:strict}.
\end{theorem}

\begin{proof}
From Theorem~\ref{thm:strict}, we have
\begin{equation}\label{eq:40}
    \norm{v^{k+1} - v^*}_H^2
    \le \norm{v^k - v^*}_H^2 - \norm{v^k - v^{k+1}}_H^2.
\end{equation}
According to Theorem~\ref{thm:monotone}, the sequence $\{\norm{v^k - v^*}_H^2\}$
is monotonically non-increasing; therefore,
\begin{equation}\label{eq:41}
    \norm{v^{t+1} - v^*}_H^2
    \le \norm{v^k - v^*}_H^2
    \le \norm{v^0 - v^*}_H^2,\quad \forall k \le t.
\end{equation}
Summing~\eqref{eq:40} for $k = 0, 1, \ldots, t$ and using~\eqref{eq:41},
\[
\sum_{k=0}^t \norm{v^k - v^{k+1}}_H^2
\le \norm{v^0 - v^*}_H^2 - \norm{v^{t+1} - v^*}_H^2
\le \norm{v^0 - v^*}_H^2.
\]
Since $\{\norm{v^k - v^{k+1}}_H^2\}$ is monotonically non-increasing,
\[
(t+1)\norm{v^t - v^{t+1}}_H^2
\le \sum_{k=0}^t \norm{v^k - v^{k+1}}_H^2
\le \norm{v^0 - v^*}_H^2,
\]
which yields~\eqref{eq:39}.
\end{proof}

The conclusion shows that, in the non-ergodic sense, the worst-case convergence rate of
ADMM for cosparse signal reconstruction is $O(1/t)$, where $t$ is the number of
iterations. This result demonstrates that the current iterate $v^k$ converges to the
solution set at the same asymptotic rate as the ergodic average $\bar{v}^t$, but provides
a stronger guarantee since it applies to the actual algorithmic output rather than a
constructed average. The key enabling factor is the monotonic non-increasing property
(Theorem~\ref{thm:monotone}), which distinguishes the non-ergodic analysis from its
ergodic counterpart.

\section{Numerical Experiments}
\label{sec:experiments}

This section presents numerical experiments to verify the ADMM method for cosparse
signal reconstruction and validate the non-ergodic convergence theory of
Section~\ref{sec:nonergodic}. All experiments were run on a HP Pro Tower 200 G9 Desktop PC with
13th Gen Intel(R) Core(TM) i5-13500 CPU (3.6 GHz), 8 GB RAM, Windows 11, and MATLAB R2025b.
The ADMM algorithm for solving~\eqref{eq:4} is summarized in Algorithm~1.

\begin{algorithm}[htbp]
\caption{ADMM for Cosparse Signal Reconstruction}
\label{alg:admm}
\begin{algorithmic}[1]
\State Task: Solve the cosparse signal reconstruction problem~\eqref{eq:4}
\State Initialization: $x^0 = 0$, $z^0 = 0$, $\lambda^0 = 0$, generate $D$, $\alpha$, $\beta$
\While{termination criterion not met}
    \State $z^{k+1} = \argmin_z L_{\beta}(z, x^k, \lambda^k)$
    \State $x^{k+1} = \argmin_x L_{\beta}(z^{k+1}, x, \lambda^k)$
    \State $\lambda^{k+1} = \lambda^k - \beta(Dx^{k+1} - z^{k+1})$
    \State $k \leftarrow k + 1$
\EndWhile
\State Check convergence conditions
\State \textbf{Output:} $z^{k+1}$, $x^{k+1}$, $\lambda^{k+1}$
\end{algorithmic}
\end{algorithm}
The maximum number of iterations was set to $\text{max\_iter} = 1000$, the penalty
parameters were set to $\alpha = 1 $ and $\beta = 1 \times 10^{-2}$,
the allowable convergence error was $\text{tol} = 1 \times 10^{-6}$, the noise level
was $\epsilon = 1 \times 10^{-2}$. $x_{\text{true}}$ represents the constructed true
cosparse signal, and $x_{\text{ADMM}}$ represents the signal reconstructed by ADMM.

In the experiments, we record the residual
$\norm{r^{k+1}}_2 = \norm{Dx^{k+1} - z^{k+1}}_2$ and the error
$\norm{x^{k+1} - x^k}_2$ at each iteration. The convergence curves confirm that the
theoretical non-ergodic rate $O(1/t)$ holds in practice.
The experimental results for two sets of parameters are shown in Figs.~1--6.

\begin{figure}[H]
    \centering
    \includegraphics[width=1\textwidth]{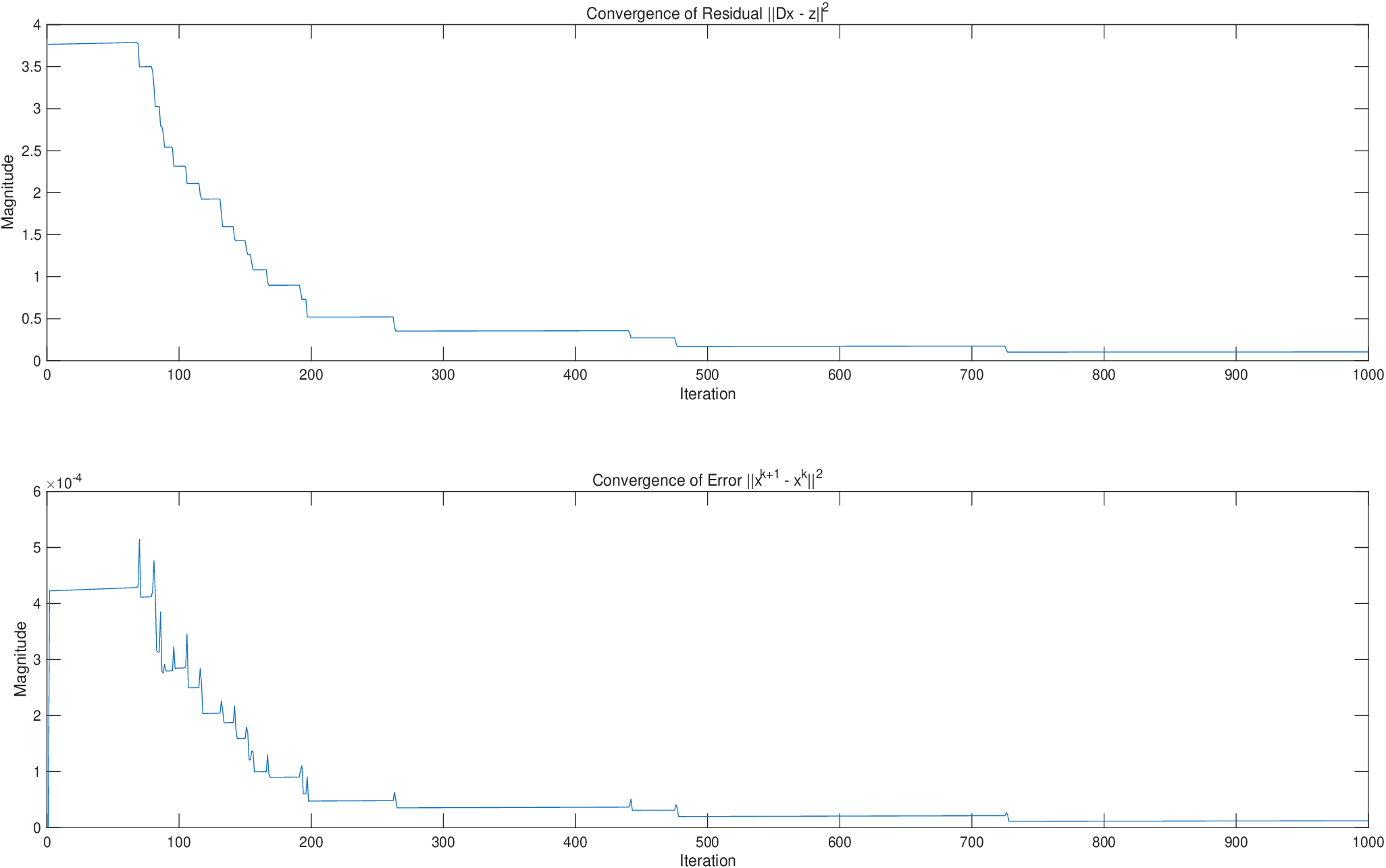}
    \caption{Residual error and convergence curve where $d = 128$, $n = 140$,
          $m = 256$, $\ell = 118$.}
    \label{fig1}
\end{figure}

\begin{figure}[H]
    \centering
    \includegraphics[width=1\textwidth]{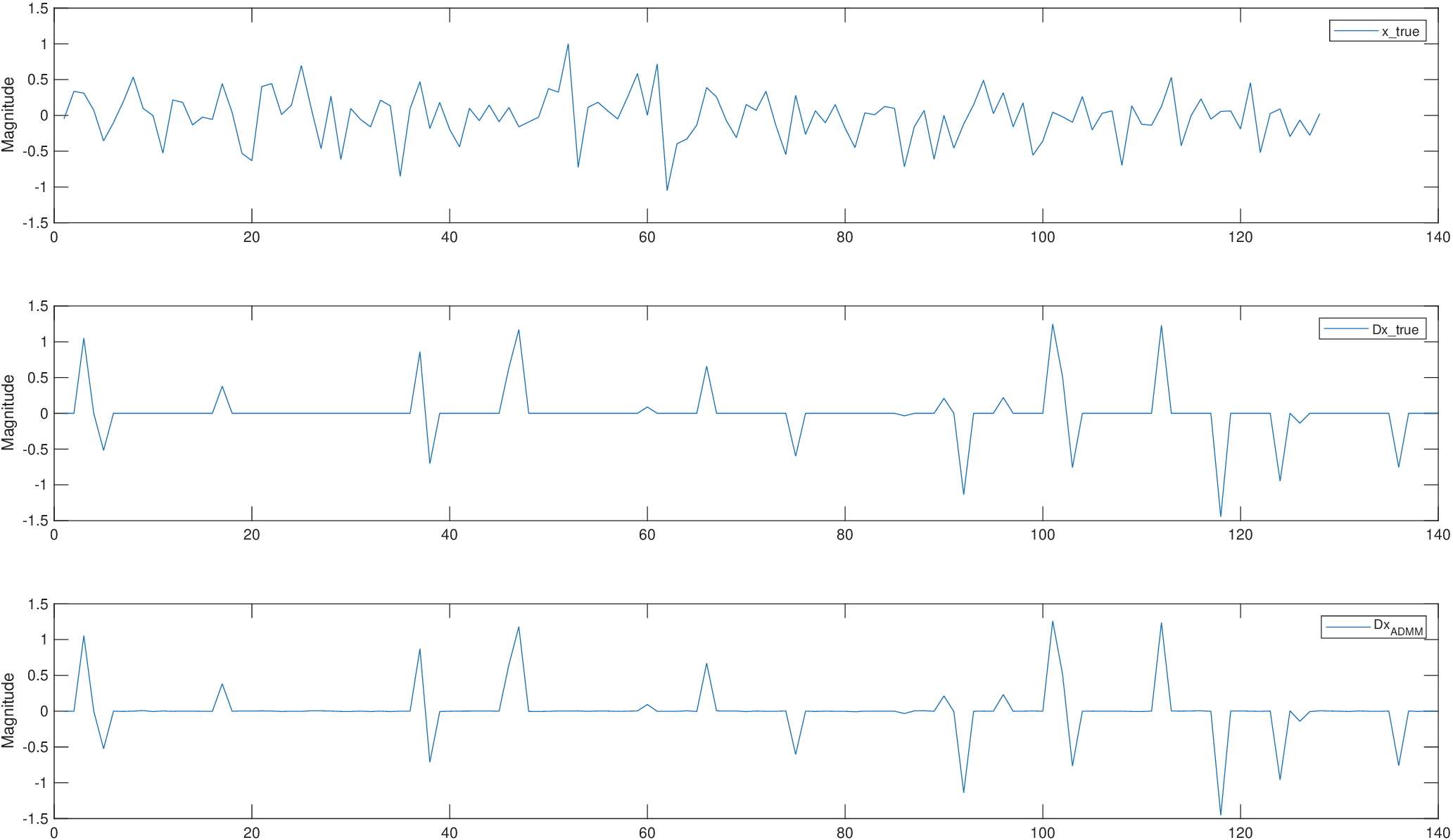}
    \caption{ Cosparse signal analysis results where $d = 128$, $n = 140$,
          $m = 256$, $\ell = 118$.}
    \label{fig2}
\end{figure}

\begin{figure}[H]
    \centering
    \includegraphics[width=1\textwidth]{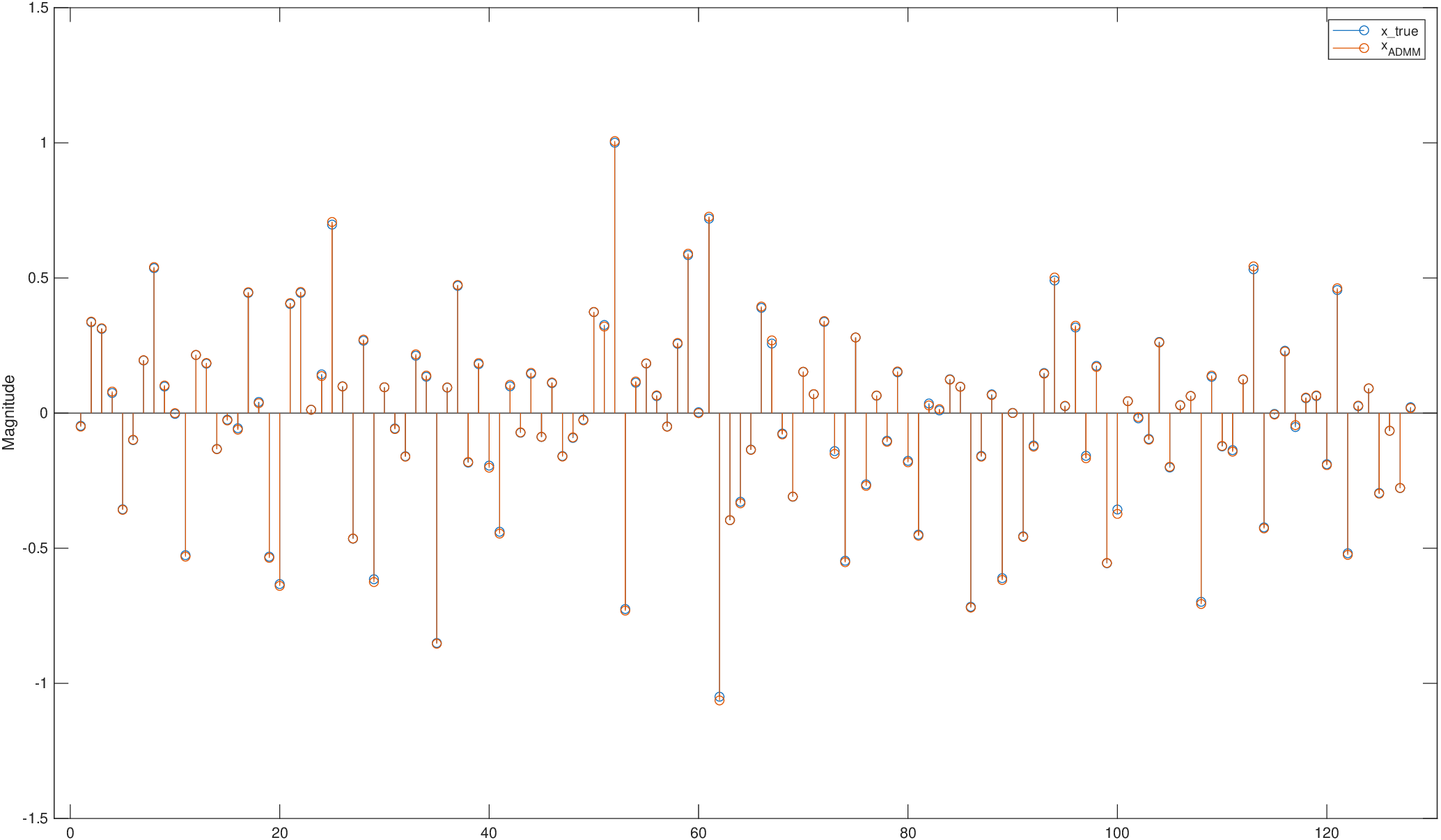}
    \caption{Reconstruction results for cosparse signals with $d = 128$, $n = 140$,
          $m = 256$, $\ell = 118$.}
    \label{fig3}
\end{figure}

\begin{figure}[H]
    \centering
    \includegraphics[width=1\textwidth]{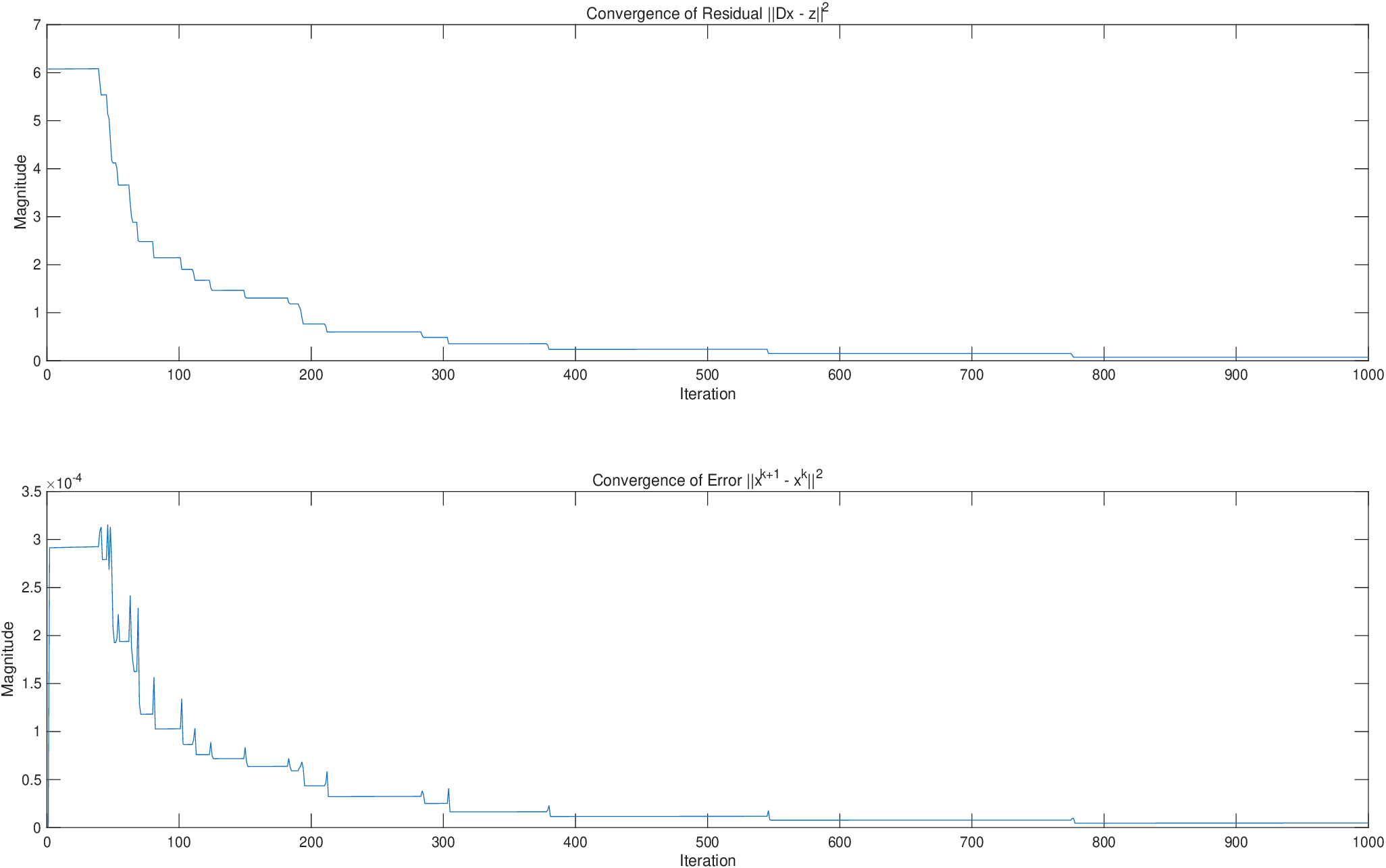}
    \caption{Cosparse signal analysis results where $d = 256$, $n = 260$,
          $m = 512$, $\ell = 236$.}
    \label{fig4}
\end{figure}

\begin{figure}[H]
    \centering
    \includegraphics[width=1\textwidth]{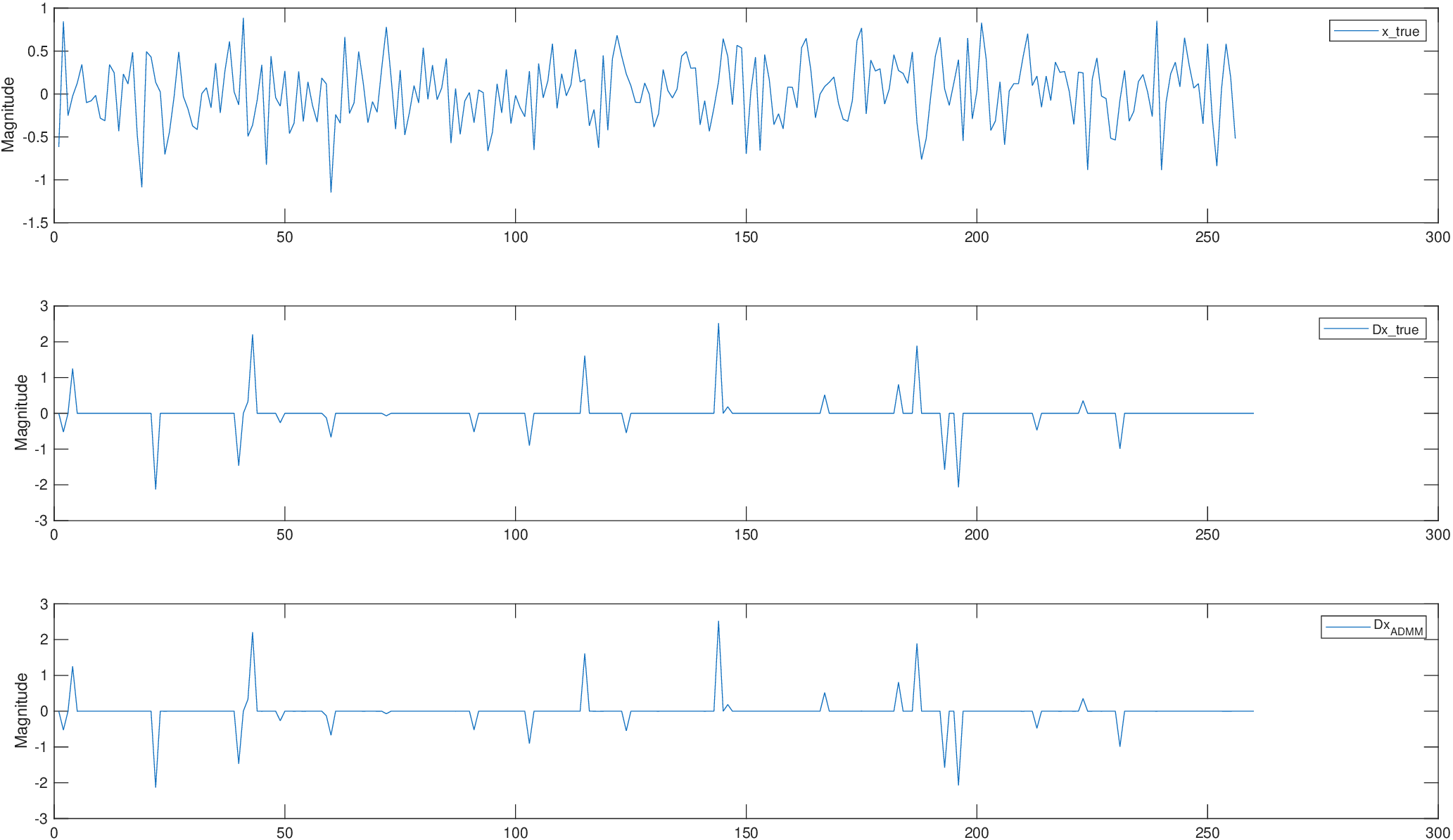}
    \caption{ Cosparse signal analysis results where $d = 256$, $n = 260$,
          $m = 512$, $\ell = 236$.}
    \label{fig5}
\end{figure}

\begin{figure}[H]
    \centering
    \includegraphics[width=1\textwidth]{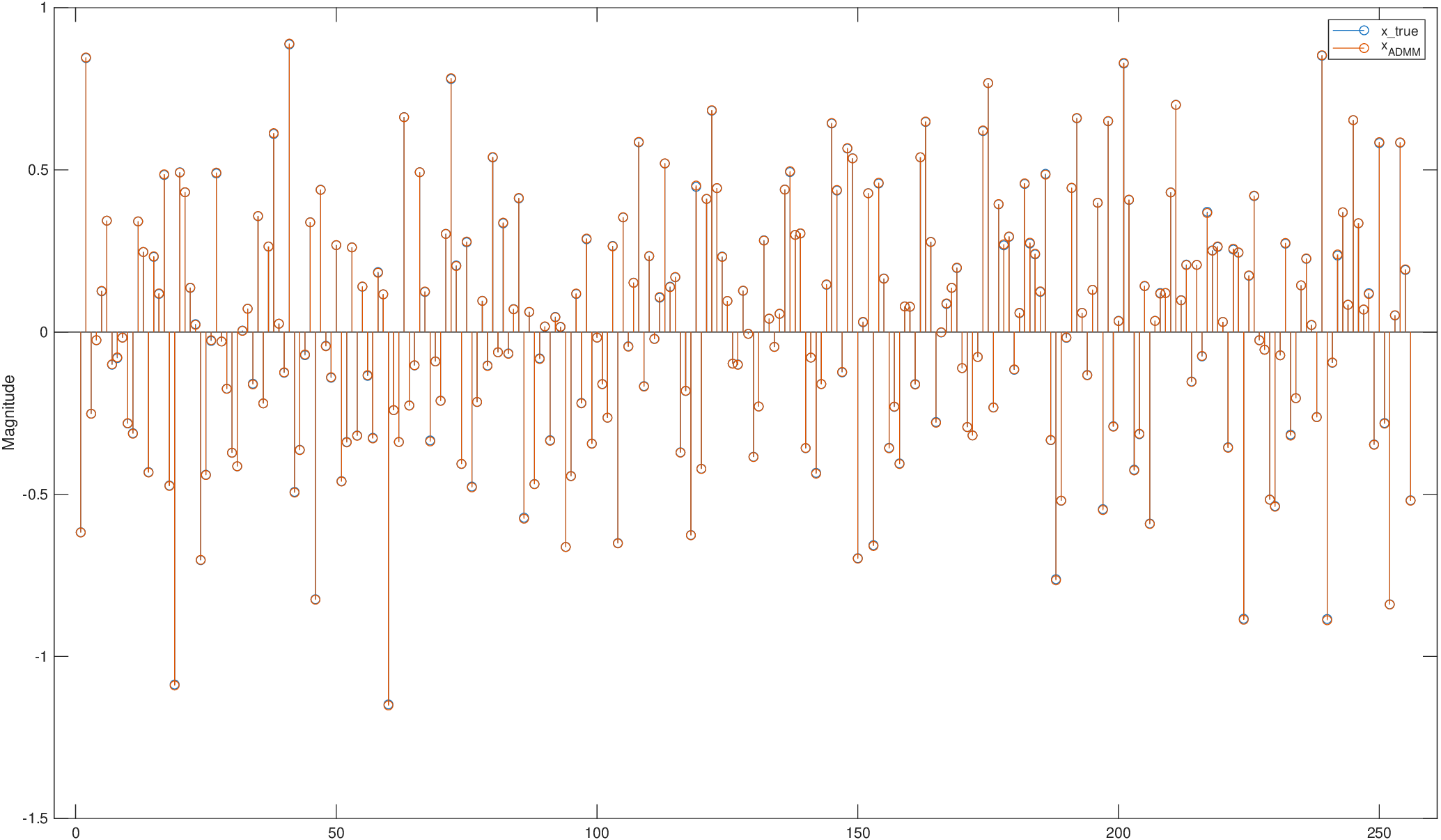}
    \caption{ Reconstruction results for cosparse signals with $d = 256$, $n = 260$,
          $m = 512$, $\ell = 236$.}
    \label{fig6}
\end{figure}

Figures 1–6 illustrate the numerical performance of ADMM for cosparse signal reconstruction under the parameter setting  $d=128$, $n=140$, $m=256$ and $\ell=118$ and $d=256$, $n=260$, $m=512$ and $\ell=236$.

Figure 1 and Figure 4 shows the convergence behavior of the algorithm. The residual $\|Dx^{k}-z^{k}\|_{2}^{2}$ decreases monotonically and approaches zero within approximately 600–700 iterations. Simultaneously, the successive error $\|x^{k+1}-x^{k}\|_{2}^{2}$ exhibits a rapid decay and stabilizes at a level on the order of $10^{-5}$ after roughly 400 iterations. These curves confirm that the iterative sequence generated by ADMM converges in practice, which is consistent with the non-ergodic $O(1/t)$ rate established in the theoretical analysis.

Figure 2 and Figure 5 displays the true cosparse signal $x_{\mathrm{true}}$ together with its analysis coefficients $Dx_{\mathrm{true}}$ and the coefficients recovered by ADMM. The coefficient vector $Dx_{\mathrm{true}}$ exhibits a clear cosparse structure, with the majority of entries close to zero and only a limited number of significant nonzero components. The recovered coefficients $Dx_{\mathrm{ADMM}}$ closely match $Dx_{\mathrm{true}}$ in both the locations and the amplitudes of the nonzero entries, indicating that the cosparse structure is accurately preserved by the algorithm.

Figure 3 and Figure 6 further compares the true signal $x_{\mathrm{true}}$ with the reconstructed signal $x_{\mathrm{ADMM}}$ in the time domain. The two signals nearly coincide at almost all sample points, demonstrating that ADMM achieves high reconstruction accuracy.

Overall, the numerical results verify both the practical effectiveness of ADMM for cosparse signal reconstruction and the correctness of the non-ergodic convergence theory developed in this paper.

\section{Conclusions}
\label{sec:conclusions}

This paper investigates the non-ergodic convergence rate of ADMM for cosparse signal
reconstruction within the VI framework. While previous work~\cite{liu2024} established
an $O(1/t)$ rate in the ergodic sense, this paper provides a stronger guarantee: the
current iterate of ADMM converges to the optimal solution set at the same asymptotic
rate in the non-ergodic sense. The main results include:
(i) the strict contraction property;
(ii) the monotonic non-increasing property of $\norm{v^k - v^*}_H^2$
(Theorem~\ref{thm:monotone}); and
(iii) the worst-case $O(1/t)$ non-ergodic convergence rate
(Theorem~\ref{thm:nonergodic}).
The analysis is independent of RIP and $D$-RIP conditions, relying only on $D^T D = I$.
Numerical experiments confirm the theoretical findings.
In future work, we plan to explore symmetric ADMM
variants~\cite{jia2021,bai2022,bai2021,bai2018,alves2020} and accelerated convergence
through inertial or adaptive penalty techniques for cosparse signal reconstruction.


\backmatter

\section*{Statements and Declarations}

\bmhead{Funding}
This work was supported by the National Natural Science Foundation of China under grant No. 62103136 and No. 12226323, Henan Province Undergraduate College Youth Backbone Teacher Training Program under grant No. 2025GGJS061, Henan Provincial Key Scientific Research Project of Higher Education Institutions, grant No. 25A120007 and 25B120002, and Henan Provincial Science and Technology Research Project under grant No. 242102210129.

\bmhead{Competing interests}
The authors declare that they have no known competing financial interests or personal relationships that could have appeared to influence the work reported in this paper.



\end{document}